\documentclass[12pt]{article}
\usepackage{amssymb}

\usepackage{graphicx}
\usepackage{graphicx}

\begin{document}

\author{Francis OGER}
\title{\textbf{Equivalence \'{e}l\'{e}mentaire entre pavages}}
\date{}

\begin{center}
\textbf{The number of paperfolding curves}

\textbf{in a covering of the plane}

\bigskip

Francis OGER
\end{center}

\bigskip

\bigskip

\noindent ABSTRACT. These results complete our previous paper in the same
Journal (vol. 42, pp. 37-75). Let $\mathcal{C}$ be a covering of the plane
by disjoint complete folding curves which satisfies the local isomorphism
property. We show that $\mathcal{C}$ is locally isomorphic to an essentially
unique covering generated by an $\infty $-folding curve. We prove that $%
\mathcal{C}$ necessarily consists of $1$, $2$, $3$, $4$ or $6$ curves. We
give examples for each case; the last one is realized if and only if $%
\mathcal{C}$ is generated by the alternating folding curve or one of its
successive antiderivatives. We also extend the results of our previous paper
to another class of paperfolding curves introduced by M. Dekking.\bigskip


\bigskip

\bigskip

\textbf{1. Definitions, context and main results.}

\bigskip

We use the definitions, the notations and the results of [4]. In order to
simplify some notations, we identify $%
\mathbb{R}
^{2}$ with $\mathbb{C}$ and $%
\mathbb{Z}
^{2}$\ with the set $%
\mathbb{Z}
+i%
\mathbb{Z}
$\ of Gaussian integers.

We consider the sequences $(a_{k})_{1\leq k\leq n-1}$, $(a_{k})_{k\in 
\mathbb{N}
^{\ast }}$, $(a_{k})_{k\in 
\mathbb{Z}
}$ with $a_{k}=\pm 1$ for each $k$ and the associated \emph{curves} $%
(C_{k})_{1\leq k\leq n}$ , $(C_{k})_{k\in 
\mathbb{N}
^{\ast }}$, $(C_{k})_{k\in 
\mathbb{Z}
}$ such that:

\noindent a) each \emph{segment} $C_{k}$ is an oriented interval $%
[x_{k},x_{k}+\varepsilon _{k}]$\ with $x_{k}\in 
\mathbb{Z}
+i%
\mathbb{Z}
$\ and\ $\varepsilon _{k}\in \left\{ 1,-1,i,-i\right\} $;

\noindent b) if $C_{k}$ and $C_{k+1}$ exist, then $x_{k}+\varepsilon
_{k}=x_{k+1}$ and $\varepsilon _{k+1}=i^{a_{k}}\varepsilon _{k}$; moreover
the curve is ``rounded'' in $x_{k+1}$\ so that it does not pass through that
point.

\noindent We call \emph{complete curves} the curves $(C_{k})_{k\in 
\mathbb{Z}
}$.

We say that a set of curves $\mathcal{C}$ \emph{covers} the plane\ (resp.
the square $S=\{x+u+iv\mid u,v\in \lbrack 0,m]\}$\ with $x\in 
\mathbb{Z}
+i%
\mathbb{Z}
$\ and\ $m\in 
\mathbb{N}
^{\ast }$) if each nonoriented interval $[z,z+\varepsilon ]\subset \mathbb{C}
$\ (resp. $[z,z+\varepsilon ]\subset S$) with $z\in 
\mathbb{Z}
+i%
\mathbb{Z}
$\ and\ $\varepsilon \in \{1,i\}$ is the support of exactly one segment of
one curve of $\mathcal{C}$. We call a \emph{covering} any set of complete
curves which covers the plane.

For each $n\in 
\mathbb{N}
^{\ast }$,\ we consider the $n$\emph{-folding sequences} $(a_{k})_{1\leq
k\leq 2^{n}-1}$ and the associated $n$\emph{-folding curves} $(C_{k})_{1\leq
k\leq 2^{n}}$, obtained by folding $n$ times a strip of paper in two, each
time possibly to the left or to the right, then unfolding it with right
angles. These curves, rounded as it is mentioned above, are self-avoiding.

We also consider the $\infty $\emph{-folding sequences} $(a_{k})_{k\in 
\mathbb{N}
^{\ast }}$, where each $(a_{1},\ldots ,a_{2^{n}-1})$ is a $n$-folding
sequence, and the associated $\infty $\emph{-folding curves}.

These two types of folding sequences and curves have been considered by
various authors (see for instance [1] and [3]).

In [4], we introduced the \emph{complete folding sequences} $(a_{k})_{k\in 
\mathbb{Z}
}$, where each $(a_{k+1},\ldots ,a_{k+l})$ is a subsequence of a $n$-folding
sequence for an integer $n$, and the associated \emph{complete folding curves%
}. For each $\infty $-folding sequence $S$, $(\overline{S},+1,S)$ and $(%
\overline{S},-1,S)$\ are complete folding sequences.

One motivation for introducing them came from two \emph{plane filling}
properties which are mentioned by various authors:

First, by [4, Th. 3.1], for each $m\in 
\mathbb{N}
^{\ast }$, there exists $n\in 
\mathbb{N}
^{\ast }$ such that, for each integer $p\geq n$, each $p$-folding curve
covers a square $\left[ x,x+m\right] \times \left[ y,y+m\right] $. It
follows that each $\infty $-folding or complete folding curve covers
arbitrarily large squares.

Second, by [4, Th. 3.15], for each $\infty $-folding curve $C$ associated to
an $\infty $-folding sequence $S$, the $4$ curves obtained from $C$ by
rotations of angles $0$, $\pi /2$, $\pi $, $3\pi /2$ around its origin are
disjoint. They can be connected in two different ways in order to form $2$
complete folding curves, both associated to $(\overline{S},+1,S)$ or both
associated to $(\overline{S},-1,S)$.

These $2$ curves form a covering, except for a particular class of $\infty $%
-folding curves. In that case, by [4, Th. 3.15], they can be completed in a
unique way to form a covering by $6$ folding curves. The $4$ other curves
are not associated to sequences of the form $(\overline{T},+1,T)$ or $(%
\overline{T},-1,T)$.

All these coverings satisfy the local isomorphism property. According to
Theorem 4 below, any covering by folding curves which satisfies the local
isomorphism property is locally isomorphic to one of them.

More generally, it follows from [4, Th. 3.10] that each complete folding
curve\ can be completed in an essentially unique way into a covering by
folding curves which satisfies the local isomorphism property. By Theorem 3
below, in such a covering, each pattern, i.e. each bounded set of curves,
appears with a well determined density.

By [4, Th. 3.12], any such covering contains at most $6$ curves. In the last
part of Section 2, we show that it can actually contain $1$, $2$, $3$ or $4$
curves. We also prove that it cannot contain $5$ curves and that it contains 
$6$ curves only in the particular case described above. We also consider the
following question: If $\mathcal{C}$\ is such a covering, for which
intergers $n$ does there exist a covering by $n$ curves which is locally
isomorphic to $\mathcal{C}$?

In [2], M. Dekking considers another notion of folding sequence. He calls a 
\emph{folding string} any sequence $(a_{1},\ldots ,a_{m-1})$ in $\left\{
-1,+1\right\} $. He introduces the \emph{folding convolution} $S\ast T$ of
two folding strings $S,T$. For each folding string $S$, he defines the
sequences $S^{\ast n}$\ with $S^{\ast 1}=S$\ and\ $S^{\ast (n+1)}=S^{\ast
n}\ast S$\ for each $n\in 
\mathbb{N}
^{\ast }$.\ Then he considers $S^{\ast \infty }=\cup _{n\in 
\mathbb{N}
^{\ast }}S^{\ast n}$.

He gives necessary and sufficient conditions on $S$ for the following
properties:

\noindent 1) $S^{\ast \infty }$\ is self-avoiding;

\noindent 2) each curve associated to $S^{\ast \infty }$\ covers arbitrarily
large squares;

\noindent 3) $S$ is \emph{perfect} in the sense that, for each curve $C$\
associated to $S^{\ast \infty }$, the plane is covered by the $4$ curves
obtained from $C$ by rotations of angles $0$, $\pi /2$, $\pi $, $3\pi /2$
around its origin.

\noindent More precisely, each of these properties is characterized by such
a condition among the sequences $S$ which satisfy the previous ones.

Many examples of $\infty $-folding sequences are actually constructed in
that way, including the positive folding sequence associated to the dragon
curve (see [4, Example 3.13]) and the alternating folding sequence (see [4,
Example 3.14]). Some of them are used in the present paper.

Some results similar to those of [4] are also true for folding sequences and
curves in Dekking's sense:

For each folding string $S=$ $(a_{1},\ldots ,a_{m-1})$, we call a \emph{%
complete }$S$\emph{-folding sequence} any sequence $T=(b_{k})_{k\in 
\mathbb{Z}
}$ such that $S^{\ast \infty }$ contains a copy of each $(b_{k+1},\ldots
,b_{k+l})$. A \emph{complete }$S$\emph{-folding curve} is a curve associated
to such a sequence.

It follows from the properties of $\ast $ that $(\overline{S}^{\ast \infty
},+1,S^{\ast \infty })$\ and $(\overline{S}^{\ast \infty },-1,S^{\ast \infty
})$\ are complete $S$-folding sequences,\ that each complete $S$-folding
sequence satisfies the local isomorphism property and that any two such
sequences are locally isomorphic.

The following results are true for each folding string $S=$ $(a_{1},\ldots
,a_{m-1})$\ which satisfies the properties 1), 2) above. Their proofs are
similar to [4]:

Any complete $S$-folding curve $C$ is self-avoiding and covers arbitrarily
large squares. We have a $\emph{derivation}$\ on $C$\ such that $m$
consecutive segments of $C$ are replaced with one segment. The derivative $%
C^{\prime }$ of $C$ is also a complete $S$-folding curve.

Each complete $S$-folding curve can be completed, in an essentially unique
way, into a covering by such curves which satisfies the local isomorphism
property. Moreover, all the coverings obtained in that way from complete $S$%
-folding curves are locally isomorphic. It would be interesting to determine
the number of curves which can appear in such a covering.

The $4$ curves considered in the property 3) are disjoint. They can be
connected in two different ways in order to form $2$ complete $S$-folding
curves associated to $(\overline{S}^{\ast \infty },+1,S^{\ast \infty })$, or 
$2$ complete $S$-folding curves associated to $(\overline{S}^{\ast \infty
},-1,S^{\ast \infty })$. If $S$\ satisfies 3), then these $2$ curves form a
covering\ which satisfies the local isomorphism property.

If $S$\ does not satisfy 3), then the $2$ curves do not form a covering, but
they can be completed in a unique way into a covering which satisfies the
local isomorphism property. It follows from Theorem 2 below that this
covering contains exactly $6$\ curves.

A simple example of that situation is obtained with $S=(+1,-1,-1)$. Then $%
S^{\ast \infty }$\ is the alternating folding sequence. It follows from [4,
Th. 3.15] that all coverings by $6$ folding curves are obtained from $\infty 
$-folding sequences $R\ast S^{\ast \infty }$\ with $R$\ finite.

Many other examples exist for folding curves in the sense of Dekking. One of
them is given in [2, Fig. 18] with $S=(+1,-1,+1,+1,-1,-1,-1,+1,-1)$.\bigskip

\textbf{2. Detailed results and proofs.}\bigskip

First we introduce some notions which will be useful both for classical
folding curves and for folding curves in the sense of Dekking.

For each curve $C$, we denote by $\alpha (C)$\ the initial point of the
first segment of $C$\ and $\beta (C)$\ the terminal point of the last
segment of $C$, if they exist.

We consider sets $\mathcal{C}$ of disjoint self-avoiding curves\ such that,
for each $z=x+iy\in 
\mathbb{Z}
+i%
\mathbb{Z}
$\ which is an endpoint of segments of curves of $\mathcal{C}$, one of the
two following possibilities is realized, depending on the parity of $x+y$:

\noindent a) the segments\ which have $z$\ as an endpoint are among $[z,z+1]$%
, $[z,z-1]$, $[z+i,z]$, $[z-i,z]$;

\noindent b) they are among $[z+1,z]$, $[z-1,z]$, $[z,z+i]$, $[z,z-i]$.

\noindent We note that this property is necessarily true if $\mathcal{C}$
consists of one curve.

For any such sets $\mathcal{C},\mathcal{D}$, we write $\mathcal{C}\cong 
\mathcal{D}$ if there exists a translation $\tau $\ such that $\tau (%
\mathcal{C})=\mathcal{D}$. We write $\mathcal{C}<\mathcal{D}$ if $\mathcal{C}%
\neq \mathcal{D}$, if each curve of $\mathcal{C}$ is contained in a curve of 
$\mathcal{D}$, and if any consecutive segments of a curve of $\mathcal{D}$
which belong to curves of $\mathcal{C}$ are consecutive in one of them. We
write $\mathcal{C}\ll \mathcal{D}$, if $\mathcal{C}<\mathcal{D}$ and if, for
each segment $A$ of a curve of $\mathcal{C}$, the $6$ segments which form
two squares with $A$\ belong to curves of $\mathcal{D}$.

For any coverings $\mathcal{C},\mathcal{D}$ such that each $z\in 
\mathbb{Z}
+i%
\mathbb{Z}
$\ is the initial point of two opposite segments, we say that a map $\Delta $%
\ from the set of segments of $\mathcal{C}$\ to the set of segments of $%
\mathcal{D}$\ is a \emph{derivation} if:

\noindent a) there\ exists a sequence $S=(a_{1},\ldots ,a_{m-1})\subset
\left\{ +1,-1\right\} $\ such that, for each segment $A$ of a curve of $%
\mathcal{D}$, $\Delta ^{-1}(A)$ is a subcurve of a curve of $\mathcal{C}$
associated to $S$ or to $\overline{S}$, depending on the parity of $x+y$\
where $\alpha (A)=x+iy$;

\noindent b) for any consecutive segments $A,B$ of a curve of $\mathcal{D}$, 
$\Delta ^{-1}(A)$\ and $\Delta ^{-1}(B)$\ are consecutive subcurves of a
curve of $\mathcal{C}$.

\noindent c) there exists a direct similitude $\sigma $\ such that $\sigma
(\alpha (A))=\alpha (\Delta ^{-1}(A))$ and $\sigma (\beta (A))=\beta (\Delta
^{-1}(A))$ for each segment $A$ of a curve of $\mathcal{D}$.

The composition of two derivations is a derivation. If $\mathcal{C}$ is a
covering by folding curves which satisfies the local isomorphism property,
then, by [4, Prop. 3.3], the derivation defined in [4] satisfies the
properties above with $m=2$. For each folding string $S=(a_{1},\ldots
,a_{m-1})$ which satisfies the properties 1), 2) of Section 1, we have a
derivation associated to $S$ on each covering by $S$-folding curves which
satisfies the local isomorphism property.\bigskip

\noindent \textbf{Proposition 1.} Consider a covering $\mathcal{C}$ with a
derivation $\Delta :\mathcal{C}\rightarrow \mathcal{C}$, a set of curves $%
\mathcal{F}<\mathcal{C}$\ and a translation $\tau $\ such that $\tau (%
\mathcal{F})\ll \Delta ^{-1}(\mathcal{F})$.\ For each $n\in 
\mathbb{N}
$, denote by $\tau _{n}$\ the translation such that $\tau _{n}(\Delta ^{-n}(%
\mathcal{F}))=\Delta ^{-n}(\tau (\mathcal{F}))\subset \Delta ^{-n-1}(%
\mathcal{F})$. Then the inductive limit of the sets $\Delta ^{-n}(\mathcal{F}%
)$\ relative to the translations\ $\tau _{n}$\ is a covering $\mathcal{D}$\
with the same number of curves as $\mathcal{F}$. Moreover, $\mathcal{D}$\ is
locally isomorphic to $\mathcal{C}$ if $\mathcal{C}$\ satisfies the local
isomorphism property.\bigskip

\noindent \textbf{Proof.} It suffices to prove that $\tau _{n}(\Delta ^{-n}(%
\mathcal{F}))\ll \Delta ^{-n-1}(\mathcal{F})$ for each $n\in 
\mathbb{N}
^{\ast }$. We show that, for each segment $S$ of $\Delta ^{-n}(\mathcal{F})$%
, the segments $S_{1},\ldots ,S_{6}$ which form two squares with $\tau
_{n}(S)$\ all belong to $\Delta ^{-n-1}(\mathcal{F})$.

Write $U_{0}=V_{0}=\tau (\Delta ^{n}(S))$. Consider the segments $%
U_{1},U_{2},U_{3},V_{1},V_{2},V_{3}\in \Delta ^{-1}(\mathcal{F})$ such that $%
U_{0},\ldots ,U_{4}$ (resp. $V_{0},\ldots ,V_{4}$) are consecutive segments
of a square.

Consider the closed curve $A=\Delta ^{-n}(U_{0})\cup \cdots \cup \Delta
^{-n}(U_{3})$\ and the closed region $P$ limited by $A$. Note that, for each 
$i\in \{0,\ldots ,3\}$, the last segment of $A_{i}$\ and the first segment
of $A_{i+1}$\ form a right angle directed to the exterior of $P$ (here we
identify $3+1$ and $0$).\ As $\mathcal{C}$ is a covering and $\Delta
^{-n}(U_{0}),\ldots ,\Delta ^{-n}(U_{3})$\ are subcurves of curves of $%
\mathcal{C}$, it follows that no curve of $\mathcal{C}$ can cross the
frontier $A$ of $P$. Consequently, the interior of $P$ contains no segment.

We prove in the same way that the interior of the closed region $Q$ limited
by $B=\Delta ^{-n}(V_{0})\cup \cdots \cup \Delta ^{-n}(V_{3})$ contains no
segment.\ As $S_{1},\ldots ,S_{6}$ are necessarily contained in $P\cup Q$,
it follows that they belong to $A\cup B\subset \Delta ^{-n-1}(\mathcal{F})$%
.~~$\blacksquare $\bigskip

Concerning folding curves in the sense of Dekking, we have:\bigskip

\noindent \textbf{Theorem 2.} Let $S$ be a folding string which satisfies
the properties 1), 2) of Section 1, but not the property 3). Let $\mathcal{C}
$ be a covering by $S$-folding curves which satisfies the local isomorphism
property and contains a curve $C$ associated to $(\overline{S}^{\ast \infty
},\mp 1,S^{\ast \infty })$. Then $\mathcal{C}$\ contains exactly $6$\
curves.\bigskip

\noindent \textbf{Proof.} We consider a covering $\mathcal{D}$ and a
derivation $\Delta :\mathcal{C}\rightarrow \mathcal{D}$\ associated to $S$.
We have $\Delta (C)\cong C$.\ Consequently, we can suppose $\Delta (C)=C$,
which implies $\mathcal{D}=\mathcal{C}$ since $C$ can be extended into a
unique covering which satisfies the local isomorphism property. We can also
suppose without restricting the generality that $S$ begins with $+1$ and
that the two first segments of $C$\ are $[0,1]$\ and $[1,1+i]$.

It follows from [2, Th. 5] and its proof that $\mathcal{C}$\ contains at
least $6$\ curves: Two of them including $C$ are associated to $(\overline{S}%
^{\ast \infty },\varepsilon ,S^{\ast \infty })$ with $\varepsilon =\mp 1$.
Each of these $2$ curves contains $2$ of the segments $[0,1]$, $[0,-1]$, $%
[i,0]$, $[-i,0]$. On the other hand, they do not contain the segments $%
[1+i,i]$, $[-1,-1+i]$, $[-1-i,-i]$, $[1,1-i]$, which necessarily belong to $%
4 $ other curves.

We are going to prove that $\mathcal{F}\ll \Delta ^{-2}(\mathcal{F})$\ for
the set of $6$\ disjoint curves $\mathcal{F}<\mathcal{C}$\ which consists of
the $8$\ segments $[0,1]$, $[0,-1]$, $[i,0]$, $[-i,0]$, $[1+i,i]$, $%
[-1,-1+i] $, $[-1-i,-i]$, $[1,1-i]$. Then it follows from Proposition 1
applied to the derivation $\Delta ^{2}$ that $\cup _{n\in 
\mathbb{N}
}\Delta ^{-2n}(\mathcal{F})\subset \mathcal{C}$\ is a covering by $6$\
curves, and therefore $\mathcal{C}$\ contains exactly $6$\ curves.

By symmetry, it suffices to show that the $6$ segments which form $2$\
squares with $[1+i,i]$\ belong to $\Delta ^{-2}(\mathcal{F})$. Our
hypotheses imply $[0,1]\in \Delta ^{-1}([0,1])$, $[1,1+i]\in \Delta
^{-1}([0,1])$, $[i,0]\in \Delta ^{-1}([i,0])$ and $[1+i,i]\in \Delta
^{-1}([1+i,i])$.

According to [2], there exists $z\in (%
\mathbb{Z}
+i%
\mathbb{Z}
)-\{0,1,-1,i,-i\}$\ such that $z=\beta (\Delta ^{-1}([0,1]))=\alpha (\Delta
^{-1}([1,1+i]))$, $(1+i)z=\beta (\Delta ^{-1}([1,1+i]))=\alpha (\Delta
^{-1}([1+i,i]))$ and $iz=\beta (\Delta ^{-1}([1+i,i]))=\alpha (\Delta
^{-1}([i,0]))$. It follows that $\alpha (\Delta ^{-1}([1+i,i]))\neq 1+i$\
and\ $\beta (\Delta ^{-1}([1+i,i]))\neq i$. Consequently,\ $\Delta
^{-1}([1+i,i])$\ contains $[1+2i,1+i]$\ and $[i,2i]$.

It follows that $[i,0]$, $[0,1]$, $[1,1+i]$, $[1+2i,1+i]$, $[i,2i]$\ belong
to $\Delta ^{-1}(\mathcal{F})$. Now it suffices to show that $[2i,1+2i]\in
\Delta ^{-2}([1+i,i])$.

We observe that\ $\Delta ^{-1}([1+2i,1+i])$\ and\ $\Delta ^{-1}([i,2i])$\
necessarily have a common vertex since $\Delta ^{-1}([0,1])$ and $\Delta
^{-1}([1+i,i])$ have the common vertex $1+i$. As $\mathcal{C}$ is a
covering, it follows that $[2i,1+2i]\in \Delta ^{-1}([1+2i,1+i])\cup \Delta
^{-1}([1+i,i])\cup \Delta ^{-1}([i,2i])\subset \Delta ^{-2}([1+i,i])$.~~$%
\blacksquare $\bigskip

From now on, we consider a covering $\mathcal{C}$ by folding curves which
satisfies the local isomorphism property. We do not mention the orientation
of the curves when it is not necessary.

The definition of the sets $E_{n}(\mathcal{C})$ for $n\in 
\mathbb{N}
\cup \left\{ \infty \right\} $\ and $F_{n}(\mathcal{C})$ for $n\in 
\mathbb{N}
$\ is given in [4]. Their existence follows from [4, Prop. 3.3]. When there
is no ambiguity, we write $E_{n}$\ and $F_{n}$ instead of $E_{n}(\mathcal{C}%
) $ and $F_{n}(\mathcal{C})$.

For each $n\in 
\mathbb{N}
$ and each $z\in E_{n}$, the $4$ nonoriented subcurves of curves of $%
\mathcal{C}$ with endpoint $z$ and length $2^{n}$ are all obtained from one
of them by successive rotations of center $z$ and angle $\pi /2$.

We say that $z=x+iy\in 
\mathbb{Z}
+i%
\mathbb{Z}
$\ is \emph{even} (resp. \emph{odd}) if $x+y$\ is even (resp. odd).

For each $n\in 
\mathbb{N}
$ and each $u\in 
\mathbb{Z}
+i%
\mathbb{Z}
$, we have:

\noindent $F_{2n}=\{u+2^{n}v\mid v\in 
\mathbb{Z}
+i%
\mathbb{Z}
$\ even$\}$ if $u\in F_{2n}$;

\noindent $E_{2n+1}=\{u+2^{n}v\mid v\in 
\mathbb{Z}
+i%
\mathbb{Z}
$\ even$\}$ if $u\in E_{2n+1}$;

\noindent $F_{2n+1}=\{u+2^{n+1}v\mid v\in 
\mathbb{Z}
+i%
\mathbb{Z}
\}$ if $u\in F_{2n+1}$;

\noindent $E_{2n+2}=\{u+2^{n+1}v\mid v\in 
\mathbb{Z}
+i%
\mathbb{Z}
\}$ if $u\in E_{2n+2}$.

For each $u\in 
\mathbb{Z}
+i%
\mathbb{Z}
$, the translation $\tau _{u}:v\rightarrow u+v$\ preserves (resp. inverses)
the orientation of the segments for $u$\ even (resp. odd).

For $n\in 
\mathbb{N}
$ and $u\in 
\mathbb{Z}
+i%
\mathbb{Z}
$\ even, we have $\tau _{2^{n}u}(F_{2n})=F_{2n}$. For each $v\in F_{2n}$,
the connections between the $4$ segments which have $v$ as an endpoint are
preserved by $\tau _{2^{n}u}$ if and only if $u\in 2(%
\mathbb{Z}
+i%
\mathbb{Z}
)$.

For $n\in 
\mathbb{N}
$ and $u\in 
\mathbb{Z}
+i%
\mathbb{Z}
$, we have $\tau _{2^{n+1}u}(F_{2n+1})=F_{2n+1}$. For each $v\in F_{2n+1}$,
the connections between the $4$ segments which have $v$ as an endpoint are
preserved by $\tau _{2^{n+1}u}$ if and only if $u$\ is even.

It follows that, for each $n\in 
\mathbb{N}
$, each $u\in 
\mathbb{Z}
+i%
\mathbb{Z}
$ and each set of curves $\mathcal{B}<\mathcal{C}$, we have $\tau
_{2^{n+1}u}(\mathcal{B})<\mathcal{C}$ if $\mathcal{B}$ contains no pair of
consecutive segments with a common vertex in $E_{2n+1}$.

We denote by $%
\mathbb{R}
_{+}$\ (resp. $%
\mathbb{R}
_{+}^{\ast }$) the set of positive (resp. strictly positive) integers.\ For
each bounded set of curves $\mathcal{B}<\mathcal{C}$, we say that the \emph{%
density} of $\mathcal{B}$ in $\mathcal{C}$ is $d\in 
\mathbb{R}
_{+}$\ if, for each $\varepsilon \in 
\mathbb{R}
_{+}^{\ast }$, there exists $r\in 
\mathbb{R}
_{+}^{\ast }$\ such that, for each $s\in 
\mathbb{R}
_{+}^{\ast }$ with $s\geq r$\ and each $z\in \mathbb{C}$,

\noindent $s^{2}d(1-\varepsilon )<|\{\mathcal{F}<\mathcal{C}\mid \mathcal{F}%
\cong \mathcal{B}$ and $\mathcal{F}\subset \Sigma
_{s}(z)\}|<s^{2}d(1+\varepsilon )$,

\noindent where $\Sigma _{s}(z)=\{z+x+iy\mid x,y\in \lbrack 0,s]\}$.\bigskip

\noindent \textbf{Theorem 3.} Each bounded set of curves $\mathcal{B}<%
\mathcal{C}$ has a density $d>0$ in $\mathcal{C}$.\bigskip

\noindent \textbf{Remark.} As the definition of density is local, $\mathcal{B%
}$ has the same density in any covering $\mathcal{D}$ which is locally
isomorphic to $\mathcal{C}$.\bigskip

\noindent \textbf{Proof.} First we give a lower bound for $d$. As $\mathcal{C%
}$ satisfies the local isomorphism property and $\left| E_{\infty }\right|
\leq 1$, there exists a copy of $\mathcal{B}$ with no vertex in $E_{\infty }$%
. Let $n$ be the smallest integer such that some $\mathcal{A}<\mathcal{C}$
with $\mathcal{A}\cong \mathcal{B}$ has no vertex in $E_{2n+1}$. Then we
have $\tau _{2^{n+1}z}(\mathcal{A})<\mathcal{C}$ for each $z\in 
\mathbb{Z}
+i%
\mathbb{Z}
$, and therefore $d\geq 1/2^{2n+2}$.

Now we prove that $d$ exists.\ For each $u\in \mathbb{C}$\ and each $s\in 
\mathbb{R}
_{+}^{\ast }$, we consider $E_{s}(u)=\{\mathcal{F}<\mathcal{C}\mid \mathcal{F%
}\cong \mathcal{B}$ and $\mathcal{F}\subset \Sigma _{s}(u)\}$. We show that,
for $s$ large and $u,v\in \mathbb{C}$, $\left| E_{s}(v)\right| -\left|
E_{s}(u)\right| $ is small compared to $s^{2}$.

We consider two integers $r,s$ such that $s$ is large compared to $2^{r}$
and $2^{r}$ is large compared to the size of $\mathcal{B}$. There exists $%
z\in 
\mathbb{Z}
+i%
\mathbb{Z}
$\ such that $v-u-2^{r+1}z=x+iy$ with $\sup (\left| x\right| ,\left|
y\right| )\leq 2^{r}$. We have $\left| E_{s}(u)\right| \leq \left|
E_{s}(v)\right| +m+n$ where

\noindent $m=|\{\mathcal{F}<\mathcal{C}\mid \mathcal{F}\cong \mathcal{B}$, $%
\mathcal{F}\subset \Sigma _{s}(u)$ and $\tau _{2^{r+1}z}(\mathcal{F}%
)\nsubseteq \Sigma _{s}(v)\}|$\ and

\noindent $n=|\{\mathcal{F}<\mathcal{C}\mid \mathcal{F}\cong \mathcal{B}$, $%
\mathcal{F}\subset \Sigma _{s}(u)$ and $\tau _{2^{r+1}z}(\mathcal{F})\nless 
\mathcal{C}\}|$.

The integer $m$ is small compared to $s^{2}$ since $s$ is large compared to $%
2^{r}$. The integer $n$ is small compared to $s^{2}$ because $2^{r}$ is
large compared to the size of $\mathcal{B}$, and because we have $\tau
_{2^{r+1}z}(\mathcal{F})<\mathcal{C}$ if $\mathcal{F}<\mathcal{C}$ contains
no pair of consecutive segments with a common vertex in $E_{2r+1}$.

It follows that $\left| E_{s}(u)\right| -\left| E_{s}(v)\right| $ is small
compared to $s^{2}$ if it is positive. The same result is true for $\left|
E_{s}(v)\right| -\left| E_{s}(u)\right| $.~~$\blacksquare $\bigskip

From now on, we do not use the identification of $%
\mathbb{R}
^{2}$ with $\mathbb{C}$.\bigskip

\noindent \textbf{Theorem 4.} $\mathcal{C}$ is locally isomorphic to a
covering generated by a curve associated to an $\infty $-folding
sequence.\bigskip

\noindent \textbf{Proof.} For $n\in 
\mathbb{N}
$ and $(u,v)\in E_{n}(\mathcal{C})$, we denote by $\mathcal{C}_{n}(u,v)$\
the set of curves obtained from $\mathcal{C}$ by keeping\ only the segments
contained in $[u-2^{n},u+2^{n}]\times \lbrack v-2^{n},v+2^{n}]$. For $m<n$, $%
\mathcal{C}_{m}(u,v)$\ is the restriction of $\mathcal{C}_{n}(u,v)$\ to $%
[u-2^{m},u+2^{m}]\times \lbrack v-2^{m},v+2^{m}]$.

By K\"{o}nig's Lemma, there exists a sequence $(X_{n})_{n\in 
\mathbb{N}
}\in \prod_{n\in 
\mathbb{N}
}E_{n}(\mathcal{C})$ such that, for $m<n$, the translation $X_{m}\rightarrow
X_{n}$\ induces an embedding of $\mathcal{C}_{m}(X_{m})$ in $\mathcal{C}%
_{n}(X_{n})$. As $\mathcal{C}$ satisfies the local isomorphism property, the
inductive limit of $(\mathcal{C}_{n}(X_{n}))_{n\in 
\mathbb{N}
}$\ relative to these embeddings is a covering $\mathcal{D}$ which is
locally isomorphic to $\mathcal{C}$ and satisfies the local isomorphism
property. The image $X$ of the elements $X_{n}$\ in $\mathcal{D}$ belongs to 
$E_{\infty }(\mathcal{D})$. Each of the two halves of curves of $\mathcal{D}$
which start at $X$ is associated to an $\infty $-folding sequence.~~$%
\blacksquare $\bigskip

\noindent \textbf{Remark.} According to [4, Theorem 3.10], the covering $%
\mathcal{D}$ given by Theorem 4 is essentially unique: two such coverings
only differ by a translation or/and a change in the connections at the $%
E_{\infty }$ point.\bigskip

For each $(x,y)\in 
\mathbb{Z}
^{2}$, the \emph{unit square} $[x,x+1]\times \lbrack y,y+1]$\ is essentially
contained in one of the connected components of $%
\mathbb{R}
^{2}-\mathcal{C}$, but each of its $4$\ vertices can belong to that
component or to another one. We say that two unit squares $S,T$ are \emph{%
connected} if they have exactly $1$ common vertex $X$ and\ if their centers
and $X$ belong to the same component.

For each $X\in 
\mathbb{Z}
^{2}$, we write $P(X)$ if $X\in E_{2}$\ and if each unit square with vertex $%
X$\ is connected to $2$ unit squares without vertex $X$.\bigskip

\noindent \textbf{Lemma 5.} There exists $A\in 
\mathbb{Z}
^{2}$ such that $\{X\in 
\mathbb{Z}
^{2}\mid P(X)\}=A+%
\mathbb{Z}
(2,-2)+%
\mathbb{Z}
(2,2)$.\bigskip

\noindent \textbf{Proof.} For each $X\in E_{2}$, the $4$ nonoriented
subcurves of curves of $\mathcal{C}$ with endpoint $X$ and length $4$ are
all obtained from one of them by successive rotations of center $X$ and
angle $\pi /2$. Consequently, for each of them, $X$ satisfies $P$ if and
only if the second and the third segment starting from $X$ are obtained from
the first one by turning left then right, or right then left.

It follows that each $X\in E_{2}$ satisfies $P$ if and only if $X+(2,0)$
(resp. $X+(0,2)$) does not satisfy $P$.~~$\blacksquare $\bigskip

\noindent \textbf{Notation.} We denote by $\mathrm{O}$ the point $(0,0)\in 
\mathbb{R}
^{2}$.\bigskip

\noindent \textbf{Theorem 6.} One of the two following properties is true:

\noindent 1) $\mathcal{C}$ consists of $1$, $2$, $3$ or $4$ curves;

\noindent 2) $\mathcal{C}$ consists of $6$ curves\ and $\mathcal{C}$ is
generated by a curve associated to the alternating folding sequence or to
one of its primitives.\bigskip

\noindent \textbf{Proof.} By [4, Th. 3.15], if $E_{\infty }(\mathcal{C})\neq
\emptyset $, then $\mathcal{C}$ consists of $2$ curves or the property 2)
above is true. It remains to be proved that, if $E_{\infty }(\mathcal{C}%
)=\emptyset $, then $\mathcal{C}$ consists of at most $4$ curves.

For each $X\in 
\mathbb{R}
^{2}$ and each curve $D$, we denote by $\delta (X,D)$ the minimum distance
between $X$ and a vertex of $D$.\ In the proof of [4, Th. 3.12], we saw that
there exist $k\in 
\mathbb{N}
$ and $X\in 
\mathbb{R}
^{2}$ such that $\delta (X,D)<1.16$ for each $D$\ in the $k$-th derivative $%
\mathcal{C}^{(k)}$ of $\mathcal{C}$. Moreover, $\mathcal{C}$ and $\mathcal{C}%
^{(k)}$ have the same number of curves and $E_{\infty }(\mathcal{C}%
)=\emptyset $ implies $E_{\infty }(\mathcal{C}^{(k)})=\emptyset $.
Consequently, we can replace $\mathcal{C}$ with $\mathcal{C}^{(k)}$, and
therefore suppose for the remainder of the proof that there exists $(x,y)\in 
\mathbb{R}
^{2}$ such that $\delta ((x,y),C)<1.16$ for each $C\in \mathcal{C}$.

Now we apply Lemma 5 to $\mathcal{C}$. There exists $A\in 
\mathbb{Z}
^{2}$ such that $\{B\in 
\mathbb{Z}
^{2}\mid P(B)\}=A+%
\mathbb{Z}
(2,-2)+%
\mathbb{Z}
(2,2)$, and therefore $(c,d)\in 
\mathbb{Z}
^{2}$ such that $P(c,d)$ and $\left| x-c\right| +\left| y-d\right| \leq 2$.

For each $n\in 
\mathbb{N}
$, we consider the images $(c_{n},d_{n})$\ and $(x_{n},y_{n})$\ of $(c,d)$\
and $(x,y)$\ in $\mathcal{C}^{(n)}$. We have $\left| x_{n}-c_{n}\right|
+\left| y_{n}-d_{n}\right| \leq 2$. In the proof of [4, Th. 3.12], we saw
that $\delta ((x,y),C)<1.16$ for each $C\in \mathcal{C}$ implies $\delta
((x_{n},y_{n}),D)<1.16$ for each $D\in \mathcal{C}^{(n)}$.

As $E_{\infty }(\mathcal{C})=\emptyset $, there exists a maximal integer $k$%
\ such that $(c_{k},d_{k})$ satisfies $P$ for$\ \mathcal{C}^{(k)}$.\
Replacing $\mathcal{C}$ with$\ \mathcal{C}^{(k)}$ if necessary, we can
assume $k=0$. We can also assume $(c_{0},d_{0})=(c_{1},d_{1})=\mathrm{O}$.
Then we have two cases:

\noindent a) $\mathrm{O}\in F_{2}(\mathcal{C})$;

\noindent b) $\mathrm{O}\in E_{3}(\mathcal{C})$ and, in the derived covering 
$\mathcal{C}^{\prime }$, each unit square with vertex $\mathrm{O}$ is
connected to exactly one unit square without vertex $\mathrm{O}$.

Figures 1A and 1B represent these two cases. Whatever the case, there are $2$
possible dispositions for the subcurves of length $4$ with endpoint $\mathrm{%
O}$ of the curves of $\mathcal{C}$. We only consider one of them since the
other one is equivalent modulo a symmetry.\ Similarly, we only consider one
of the $2$ possible choices for the connections in $\mathrm{O}$. We note
that the ball $B((x,y),1.16)$ is contained in the interior of the square $%
(Z_{1},Z_{2},Z_{3},Z_{4})$\ because $(x,y)$\ belongs to the square $%
(X_{1},X_{2},X_{3},X_{4})$.

In Figure 1A, the connections in $Y_{1},Y_{3}$\ are imposed by the
connections in $\mathrm{O}$, since the property $\mathrm{O}\in F_{2}(%
\mathcal{C})$ implies $X_{1},X_{2},X_{3},X_{4}\in E_{3}(\mathcal{C})$.
Because of the existence of connections in $X_{1},X_{2},X_{3},X_{4}$, all
the subcurves represented are contained in at most $4$ curves of $\mathcal{C}
$. As no other curve of $\mathcal{C}$ can reach the vertices in $%
B((x,y),1.16)$, $\mathcal{C}$\ contains at most $4$ curves.

In Figure 1B, the connections in $X_{1},X_{2},X_{3},X_{4}$\ are imposed
since, in $\mathcal{C}^{\prime }$, each unit square with vertex $\mathrm{O}$%
\ is connected to only one unit square without vertex $\mathrm{O}$.
Consequently, $\mathcal{C}$\ contains at most $2$ curves with segments in
the interior of the square $(Y_{1},Y_{2},Y_{3},Y_{4})$. As at most $2$ other
curves of $\mathcal{C}$ can reach the vertices in $B((x,y),1.16)$, it
follows that $\mathcal{C}$\ contains at most $4$ curves.~~$\blacksquare $%
\bigskip

\begin{center}
\includegraphics[scale=0.60]{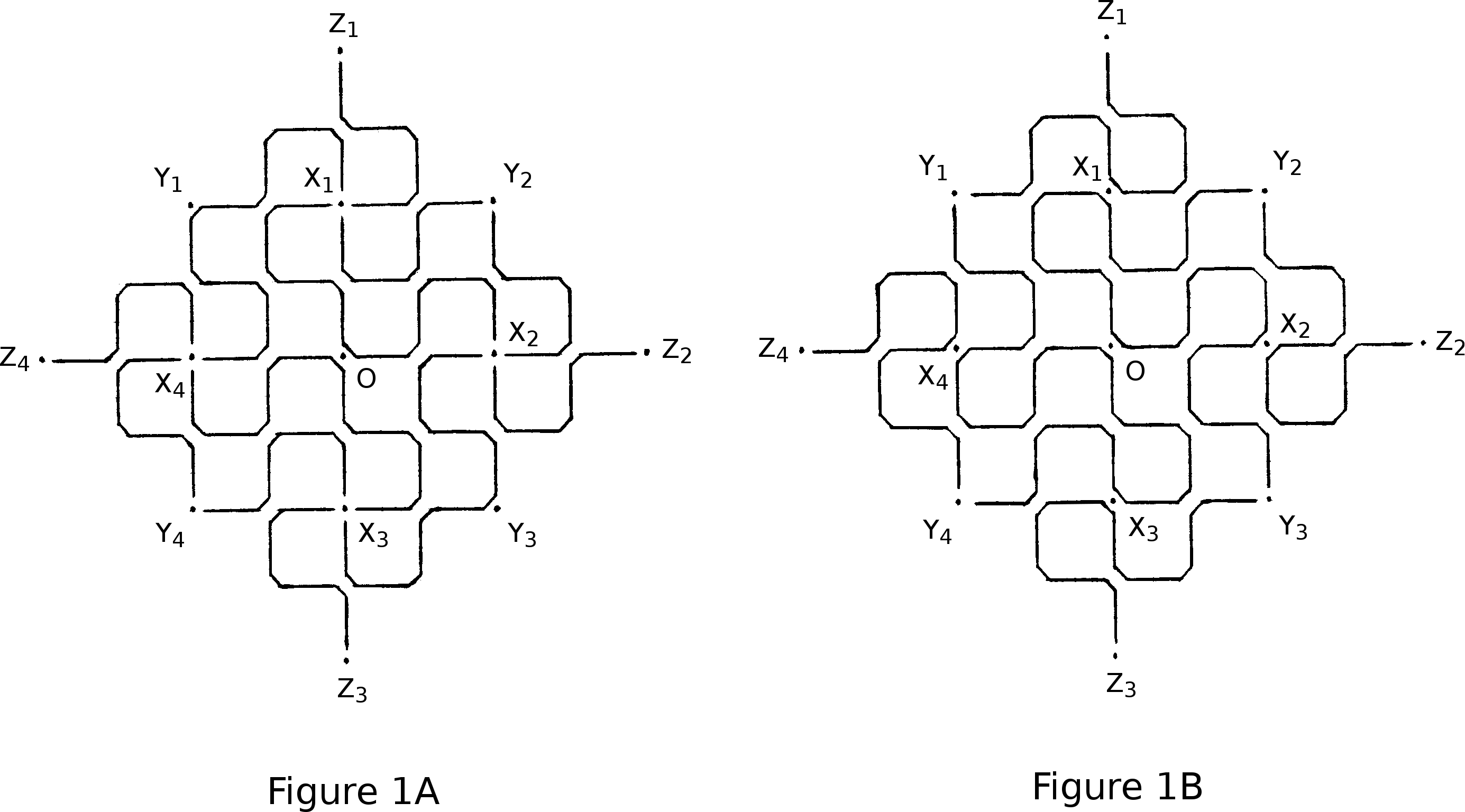}\bigskip
\end{center}

By [4, Th. 3.2, Cor. 3.6 and Th. 3.7], $\mathcal{C}$ is locally isomorphic
to a covering by $1$ curve. We also have:\bigskip

\noindent \textbf{Proposition 7.} $\mathcal{C}$ is locally isomorphic to a
covering by $2$ curves.\bigskip

\noindent \textbf{Proof.} By Theorem 4, we can suppose that $\mathcal{C}$ is
generated by a curve associated to an $\infty $-folding sequence $S$. Then,
by [4, Th. 3.15], $\mathcal{C}$ itself consists of $2$ curves except if $S$
is the alternating folding sequence or one of its primitives. So we can
suppose for the remainder of the proof that there exists $k\in 
\mathbb{N}
$ such that $S^{(k)}$ is the alternating folding sequence $T$. Then $%
\mathcal{C}^{(k)}$ is generated by a curve associated to $T$.

If there exists a covering $\mathcal{D}$ by $2$ curves which is locally
isomorphic to $\mathcal{C}^{(k)}$, then there exists a $k$-th primitive $%
\mathcal{E}$ of $\mathcal{D}$ which is locally isomorphic to $\mathcal{C}$,
and $\mathcal{E}$ also consists of $2$ curves. So we can suppose $k=0$. Then 
$\mathcal{C}$ is the covering shown in [4, Fig. 8].

For each $n\in 
\mathbb{N}
^{\ast }$ and each $r\in 
\mathbb{Z}
$, two bounded subcurves of distinct curves of $\mathcal{C}$ form a covering 
$\mathcal{C}_{n,r}$, in the sense given in the proof of [4, Example 3.8], of
the triangle $%
T_{n,r}=((0,(2r)2^{n}),(2^{n},(2r+1)2^{n}),(-2^{n},(2r+1)2^{n}))$. For each $%
s\in 
\mathbb{Z}
$, the translation $X\rightarrow X+(0,s.2^{n+1})$\ induces an isomorphism
from $\mathcal{C}_{n,r}$ to $\mathcal{C}_{n,r+s}$.

Now, for each $n\in 
\mathbb{N}
^{\ast }$, we consider the embedding $\pi _{n}:\mathcal{C}_{2n,0}\rightarrow 
\mathcal{C}_{2n,1}\subset \mathcal{C}_{2n+2,0}$\ induced by the translation $%
\tau _{n}:X\rightarrow X+(0,2^{2n+1})$. We observe that $\bigcup_{n\in 
\mathbb{N}
}\tau _{1}^{-1}(\cdots (\tau _{n}^{-1}(T_{2n+2,0}))\cdots )=%
\mathbb{R}
^{2}$. It follows that the inductive limit of\ $(\mathcal{C}_{2n,0})_{n\in 
\mathbb{N}
^{\ast }}$\ relative to the embeddings $\pi _{n}$ is a covering by $2$
curves which is locally isomorphic to $\mathcal{C}$.~~$\blacksquare $\bigskip

We do not know presently for which coverings $\mathcal{C}$ there exists a
covering $\mathcal{D}$ consisting of $3$ or $4$ curves which is locally
isomorphic to $\mathcal{C}$. However, we are going to give examples of the
two situations. Here, for any sets of curves $\mathcal{F},\mathcal{G}$, we
say that $\mathcal{F}$ is \emph{interior} to $\mathcal{G}$ if $\mathcal{F}%
\ll \mathcal{G}$.\bigskip

\noindent \textbf{Example 8.} Let $\mathcal{C}$ be a covering generated by a
dragon curve associated to the $\infty $-folding sequence $(a_{k})_{k\in 
\mathbb{N}
^{\ast }}$ with $a_{2^{n}}=+1$ for each $n\in 
\mathbb{N}
$.\ Then $\mathcal{C}$\ is locally isomorphic to two coverings by $3$
curves, one where each pair of curves has common vertices and one with two
curves separated by the third one.\bigskip

\noindent \textbf{Proof.} For the first covering, we apply Proposition 1 to
the set of curves $\mathcal{D}$ shown in Figure 2A with horizontal and
vertical segments. It is embedded in $\mathcal{C}$ since it appears in [4,
Fig. 7] between $(2,0)$ and $(2,2)$.

We consider the first primitive of $\mathcal{D}$, shown in Figures 2A and 2B
with diagonal segments, and the second primitive, shown in Figure 2B with
horizontal and vertical segments. They are also embedded in $\mathcal{C}$.
By Figure 2B, the second primitive contains the image of $\mathcal{D}$ under
a rotation by $-\pi /2$; we note that this image is not interior to it
because the condition is not satisfied for one of the segments, but it is
satisfied at the following step. Repeting this process $3$ more times, we
obtain a copy of $\mathcal{D}$ which is interior to the $8$-th primitive of $%
\mathcal{D}$.

For the second covering, we apply Proposition 1 to the set of curves $%
\mathcal{E}$ shown in Figure 3A with horizontal and vertical segments. It is
embedded in $\mathcal{C}$ since it appears in [4, Fig. 7] between $\mathrm{O}
$ and $(0,3)$.

We consider the first primitive of $\mathcal{E}$, shown in Figures 3A and
3B, and the second primitive, shown in Figure 3B. They are also embedded in $%
\mathcal{C}$. By Figure 3B, an image of $\mathcal{E}$ under a rotation by $%
\pi $ is interior to the second primitive of $\mathcal{E}$. Repeting this
process $1$ more time, we obtain a copy of $\mathcal{E}$ which is interior
to the $4$-th primitive of $\mathcal{E}$.~~$\blacksquare $\bigskip 

\begin{center}
\includegraphics[scale=0.50]{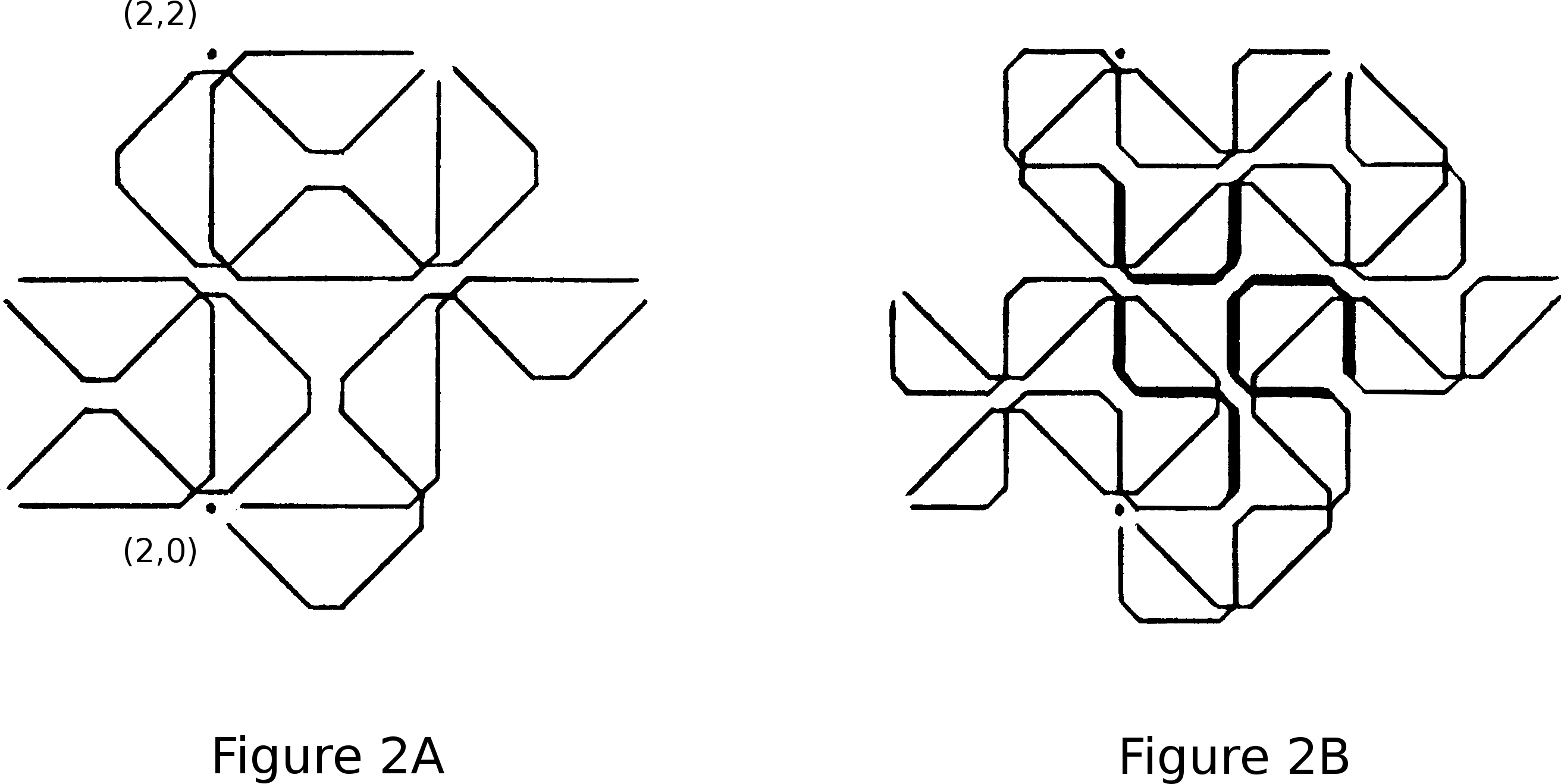}\bigskip

\includegraphics[scale=0.50]{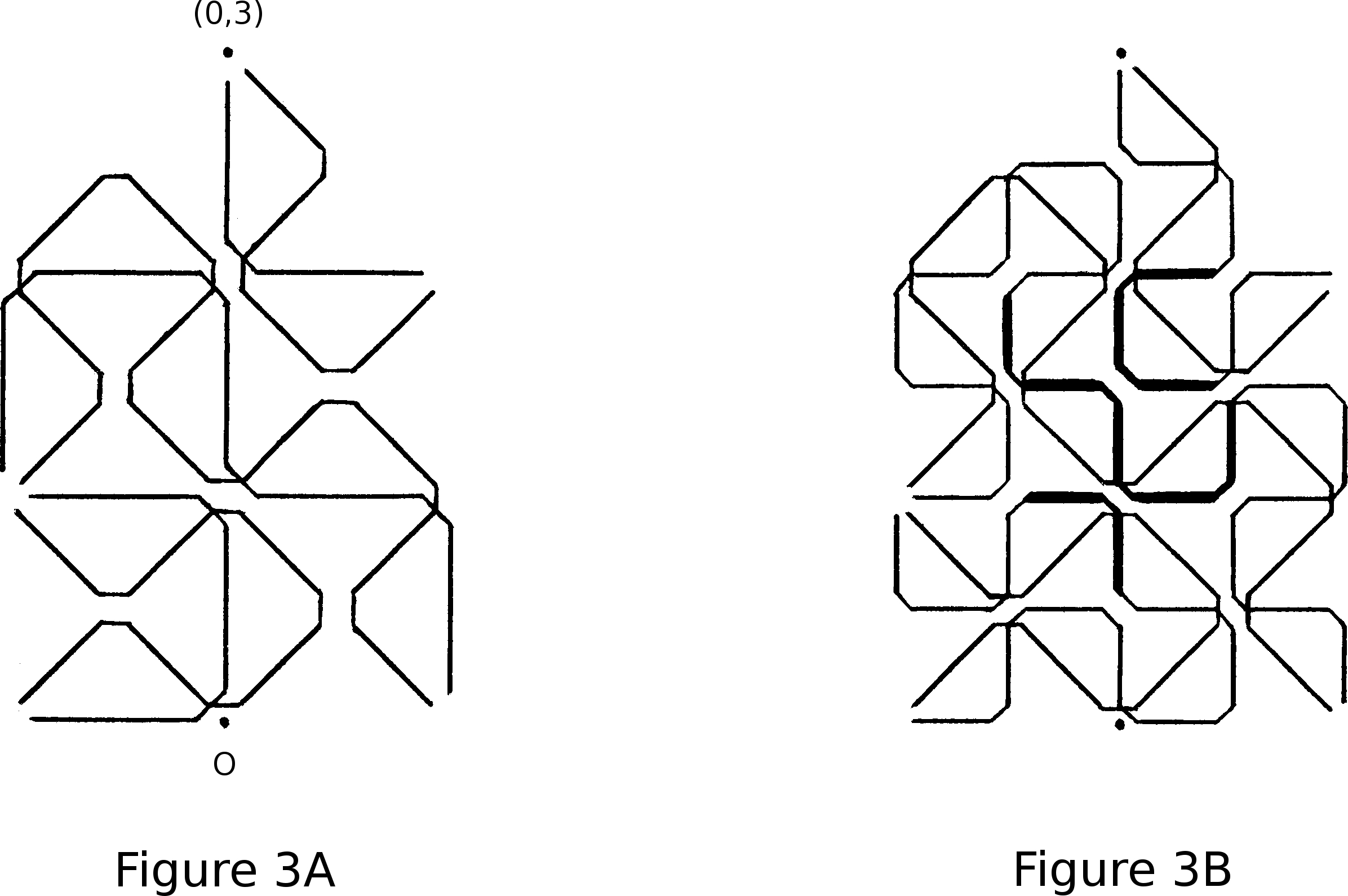}\bigskip
\end{center}

\noindent \textbf{Example 9.} Let $\mathcal{C}$ be a covering generated by a
curve associated to the $\infty $-folding sequence $(a_{k})_{k\in 
\mathbb{N}
^{\ast }}$ with $a_{2^{4n}}=a_{2^{4n+1}}=+1$\ and $%
a_{2^{4n+2}}=a_{2^{4n+3}}=-1$ for each $n\in 
\mathbb{N}
$.\ Then $\mathcal{C}$\ is locally isomorphic to a covering by $4$
curves.\bigskip

\noindent \textbf{Proof.} $\mathcal{C}$ is represented in Figure 4. We apply
Proposition 1 to the set of curves $\mathcal{D}$ shown in Figure 5A with
horizontal and vertical segments. It is embedded in $\mathcal{C}$ since it
appears in Figure 4 between $(-3,0)$ and $\mathrm{O}$.

We consider the first primitive of $\mathcal{D}$, shown in Figures 5A and
5B, and the second primitive, shown in Figure 5B. The second primitive is
embedded in the covering generated by a curve associated to the $\infty $%
-folding sequence $(a_{k})_{k\in 
\mathbb{N}
^{\ast }}$ with $a_{2^{4n}}=a_{2^{4n+1}}=-1$\ and $%
a_{2^{4n+2}}=a_{2^{4n+3}}=+1$ for each $n\in 
\mathbb{N}
$.\ According to Figure 5B, it contains a copy of the image of $\mathcal{D}$
under a reflection about the y-axis; we note that this copy is not interior
to it because the condition is not satisfied for one of the segments, but it
is satisfied at the following step.

Then, by applying a process which is the image of the previous one under a
reflection about the $y$-axis, we obtain a copy of $\mathcal{D}$ which is
interior to the $4$-th primitive of $\mathcal{D}$.~~$\blacksquare $\bigskip

\begin{center}
\includegraphics[scale=0.60]{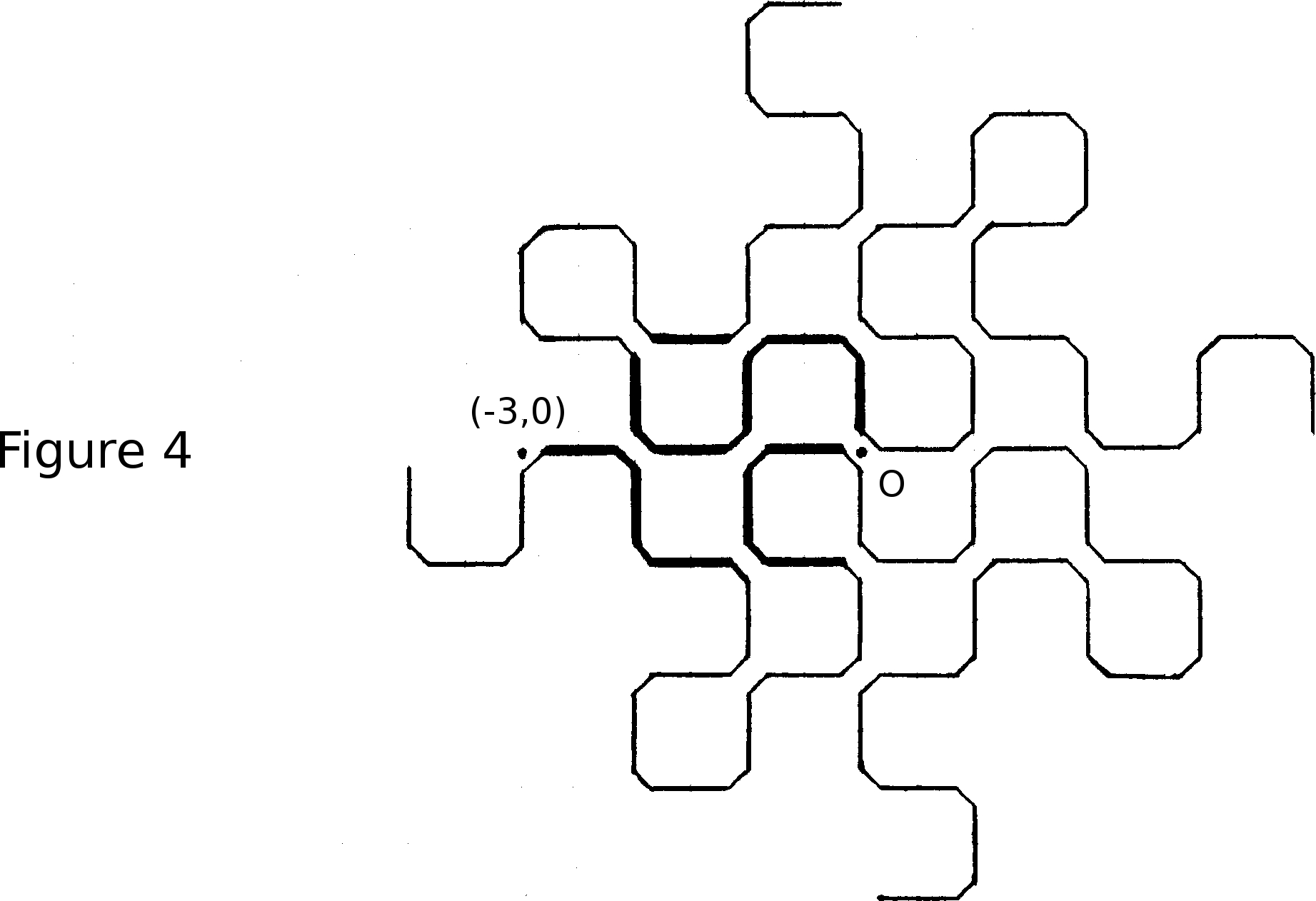}\bigskip

\includegraphics[scale=0.60]{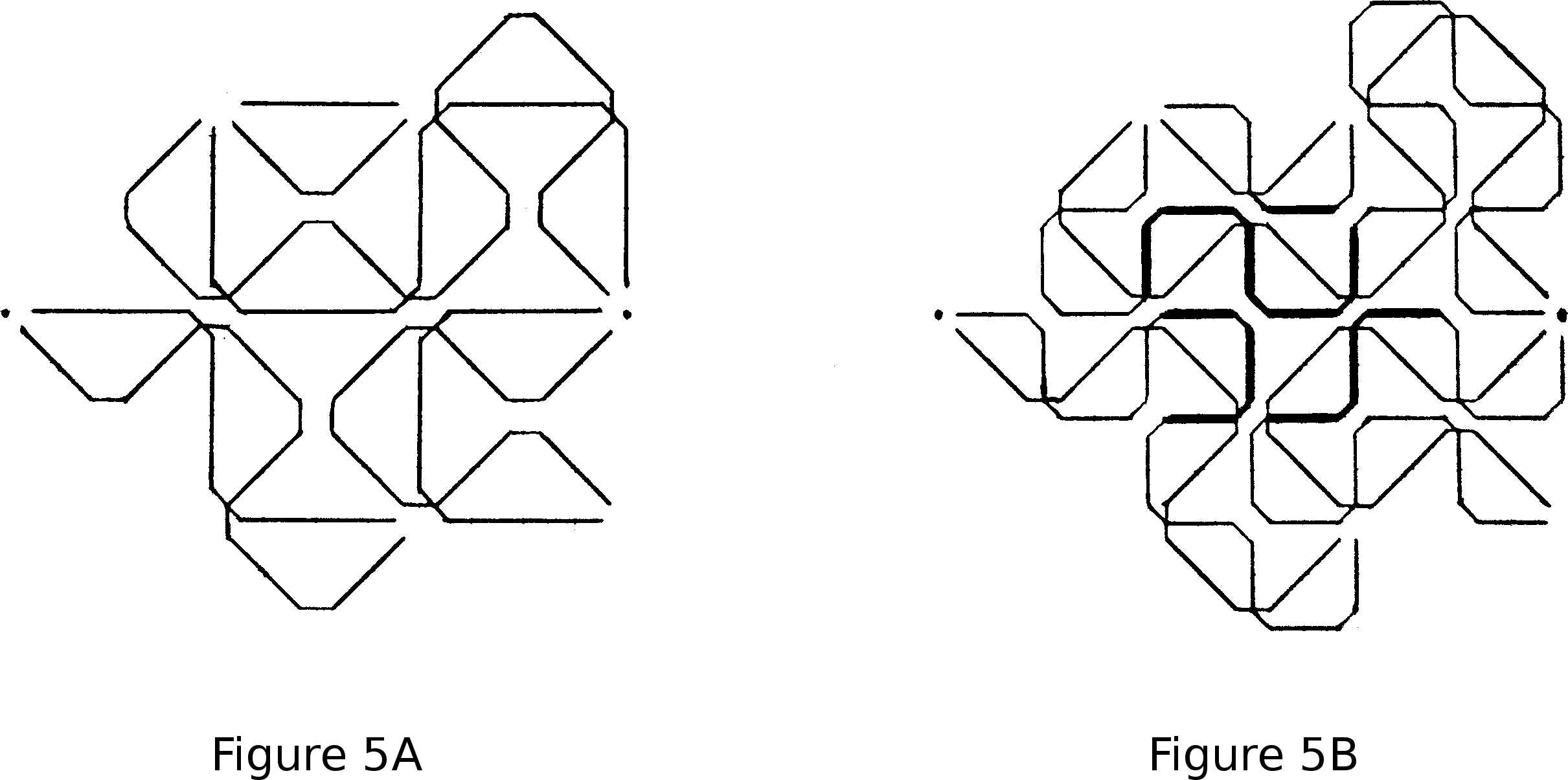}\bigskip \bigskip

\textbf{Acknowledgements}\medskip 
\end{center}

We are grateful to the referee for drawing our attention to the paper [2] by
M. Dekking.\bigskip \bigskip

\begin{center}
\textbf{References}\bigskip
\end{center}

\noindent \lbrack 1] D. Davis and D.E. Knuth, Number representations and
dragon curves I and II, J. Recreational Math. 3 (1970), 66-81 and 133-149.

\noindent \lbrack 2] M. Dekking, Paperfolding morphisms, planefilling
curves, and fractal tiles, Theoret. Comput. Sci. 414 (2012), 20-37.

\noindent \lbrack 3] M. Dekking, M. Mend\`{e}s France and A. Van Der
Poorten, Folds, Math. Intelligencer 4 (1982), 130--138.

\noindent \lbrack 4] F. Oger, Paperfolding sequences, paperfolding curves
and local isomorphism, Hiroshima Math. Journal 42 (2012), 37-75.\bigskip
\bigskip

Francis OGER

UFR de Math\'{e}matiques, Universit\'{e} Paris 7

B\^{a}timent Sophie Germain, case 7012

75205 Paris Cedex 13

France

E-mail: oger@math.univ-paris-diderot.fr.

\vfill

\end{document}